\documentclass[12pt]{article}

\usepackage[margin=1.17in]{geometry}

\usepackage{amssymb,amsmath}
\usepackage{amsthm,enumitem,mathtools,microtype}
\usepackage{hyperref}
\usepackage{mathrsfs}
\usepackage{textcomp}
\hypersetup{
	colorlinks,%
	citecolor=black,%
	filecolor=black,%
	linkcolor=black,%
	urlcolor=black
}
\usepackage[abbrev,nobysame]{amsrefs}
\usepackage{todonotes}

\usepackage{verbatim}

\usepackage{sectsty}

\usepackage[T1]{fontenc}
\usepackage{lmodern}

\usepackage{needspace}

\newtheorem{theorem}{Theorem}[section]
\newtheorem{lemma}{Lemma}[section]
\theoremstyle{remark}

\newcommand{\Sn}{\mathbb S^n}
\newcommand{\Rn}{\mathbb R^n}
\newcommand{\ds}{\,d\sigma}
\newcommand{\pc}{\frac{np}{n-p}}
\DeclareMathOperator{\supp}{supp}

\begin{document}
\title{\textbf{Sharp Sobolev inequalities under first moment constraints}\bigskip}
	
\author{\medskip  Rupert L. Frank\footnote{R.L. Frank is partially supported through the German Research Foundation grants EXC-2111-390814868, TRR 352–Project-ID 470903074 and FR 2664/3-1.},\ \ Tianling Jin\footnote{T. Jin is partially supported by NSFC grant 12122120, and Hong Kong RGC grants GRF 16304125, GRF 16303624  and GRF 16303822.}, \  \ Hua-Yang Wang\footnote{H.-Y. Wang is partially supported by NSFC grant 12601183.},\ \  Jingang Xiong\footnote{J. Xiong is partially supported by NSFC grant 12325104.}}

\date{\today}

\maketitle

\begin{abstract}
Let $n\geq2$ and $1<p<n$. We study sharp Sobolev inequalities on the round sphere $\mathbb{S}^n$ for functions $u\in W^{1,p}(\mathbb{S}^n)$ satisfying $\int_{\mathbb{S}^n}x|u|^{np/(n-p)}\,d\sigma=0$. We prove that
\[
 \|u\|_{L^{np/(n-p)}(\mathbb{S}^n)}^p \leq 2^{-p/n}S_{n,p}^p\|\nabla u\|_{L^p(\mathbb{S}^n)}^p+B\|u\|_{L^p(\mathbb{S}^n)}^p
\]
holds for all such functions with some $B=B(n,p)<\infty$ if and only if $1<p\leq2$ or $2<p<n<3p-2$, where $S_{n,p}$ is the optimal Euclidean Sobolev constant. This determines the exact range in which Aubin's inequality holds with the sharp gradient coefficient. For every $2\leq p<n$, we also establish the corresponding inequality with squared norms and gradient coefficient $2^{-2/n}S_{n,p}^2$. The gradient coefficient in each inequality is optimal. The proof combines stereographic averaging with quantitative stability of the Euclidean Sobolev inequality to obtain local second moment estimates that allow absorption of the localization errors.


\medskip

\end{abstract}

\section{Introduction}

Let $n\geq2$ and $1<p<n$. Denote by $S_{n,p}$ the optimal constant in the Euclidean Sobolev inequality
\[
 \|v\|_{L^{np/(n-p)}(\Rn)} \leq S_{n,p}\|\nabla v\|_{L^p(\Rn)}
\]
for $v\in C_c^\infty(\Rn)$. Aubin \cite{Aubin1976} and Talenti \cite{Talenti1976} determined this constant independently.

On a smooth compact Riemannian manifold $(M^n,g)$ without boundary, for every $1<p<n$ there exists $B=B(M,g,p)<\infty$ such that, for all $u\in W^{1,p}(M)$,
\[
 \begin{aligned}
 \|u\|_{L^{np/(n-p)}(M)}^p&\leq S_{n,p}^p\|\nabla_g u\|_{L^p(M)}^p+B\|u\|_{L^p(M)}^p, &&1<p\leq2,\\
 \|u\|_{L^{np/(n-p)}(M)}^2&\leq S_{n,p}^2\|\nabla_g u\|_{L^p(M)}^2+B\|u\|_{L^p(M)}^2, &&2<p<n.
 \end{aligned}
\]
Here the norms are taken with respect to the Riemannian volume measure, and the leading coefficients are optimal. The case $p=2$ was proved by Hebey--Vaugon \cite{HebeyVaugon1996}, and the cases $p\ne2$ by Aubin--Li \cite{AubinLi1999} and, independently, Druet \cite{Druet1999BestConstants}.
For $p>2$, the choice of power is essential: Druet \cite{Druet1998, Druet1999} showed that positive scalar curvature at a point can obstruct the corresponding inequality with $p$-th powers.

Let $\Sn\subset\mathbb R^{n+1}$ be the round sphere, equipped with its standard area measure $d\sigma$. We consider $u\in W^{1,p}(\Sn)$ satisfying
\begin{equation}\label{eq:constraints}
 \int_{\Sn}x_j|u(x)|^{\pc}\ds=0,\quad j=1,\cdots,n+1,
\end{equation}
where $x=(x_1,\cdots,x_{n+1})\in\Sn$. Thus the first moments of the critical density $|u|^{np/(n-p)}$ vanish. This condition excludes concentration of the entire critical mass at a single point. In his work on prescribing scalar curvature, Aubin \cite{Aubin1979} proved that, for every $\varepsilon>0$,
\begin{equation}\label{eq:aubin}
 \|u\|_{L^{np/(n-p)}(\Sn)}^p \leq\bigl(2^{-p/n}S_{n,p}^p+\varepsilon\bigr)\|\nabla u\|_{L^p(\Sn)}^p+B_\varepsilon\|u\|_{L^p(\Sn)}^p
\end{equation}
for all $u$ satisfying \eqref{eq:constraints}, where $B_\varepsilon=B(n,p,\varepsilon)<\infty$. Hang--Wang \cite{HangWang2022} gave a new proof using concentration compactness and extended the inequality to higher moment constraints. Sun--Zou \cite{SunZou2025} subsequently obtained related inequalities on compact manifolds, with an additional $\varepsilon$ in the leading coefficient.

In this paper, we study whether \eqref{eq:aubin} holds with $\varepsilon=0$ and some $B<\infty$. Two identical bubbles concentrating at antipodal points satisfy \eqref{eq:constraints} and show that $2^{-p/n}S_{n,p}^p$ cannot be decreased. The sharpness of this coefficient does not imply the existence of $B<\infty$ at $\varepsilon=0$, since \eqref{eq:aubin} gives no uniform bound for $B_\varepsilon$ as $\varepsilon\downarrow0$.

Our main theorem establishes the square-norm inequality for $2\leq p<n$ and determines the exact range in which the $p$-power inequality holds with the sharp leading coefficient.

\begin{theorem}\label{thm:main}
Let $n\geq2$ be an integer and $1<p<n$. Consider functions $u\in W^{1,p}(\Sn)$ satisfying \eqref{eq:constraints}.
\begin{enumerate}
\renewcommand{\labelenumi}{\textup{(\roman{enumi})}}
\item If $2\leq p<n$, there exists $B=B(n,p)<\infty$ such that every such $u$ satisfies
\begin{equation}\label{eq:main-square}
 \|u\|_{L^{np/(n-p)}(\Sn)}^2\leq 2^{-2/n}S_{n,p}^2\|\nabla u\|_{L^p(\Sn)}^2+B\|u\|_{L^p(\Sn)}^2.
\end{equation}

\item There exists $B=B(n,p)<\infty$ such that every such $u$ satisfies
\begin{equation}\label{eq:main-p}
 \|u\|_{L^{np/(n-p)}(\Sn)}^p\leq 2^{-p/n}S_{n,p}^p\|\nabla u\|_{L^p(\Sn)}^p+B\|u\|_{L^p(\Sn)}^p
\end{equation}
if and only if $1<p\leq2$ or $2<p<n<3p-2$.
\end{enumerate}
The coefficient of the gradient term in each inequality is optimal in the stated range.
\end{theorem}

We do not claim that the constants $B$ in these inequalities are optimal. For $p=2$ and $n\geq3$, Chang--Yang \cite{ChangYang1991} proved that, under \eqref{eq:constraints}, the gradient coefficient can be decreased below $S_{n,2}^2$ while keeping the optimal lower-order coefficient $B=|\Sn|^{-2/n}$. Li \cite{Li1996} gave a simpler proof. For $n\geq5$, Wang--Zhang \cite{WangZhang2026} obtained an explicit improved gradient coefficient, which is larger than $2^{-2/n}S_{n,2}^2$. Theorem~\ref{thm:main} gives the sharp gradient coefficient for every $n\geq3$, with some finite $B$. For $p=2$ and $u(-x)=u(x)$, Dolbeault \cite{Dolbeault2021} previously conjectured \eqref{eq:main-square} with $B=|\Sn|^{-2/n}$. Dolbeault's conjecture has recently been solved by Wei--Yan--Ye \cite{WYY}.

To prove the inequalities with the sharp gradient coefficients, we need quantitative control of the localization errors. For a concentrating sequence at the optimal energy, concentration compactness identifies two equal antipodal masses, but does not quantify the errors in the localized masses. For $p=2$, Dolbeault--Esteban--Loss \cite{DolbeaultEstebanLoss2017} related improved spherical Sobolev inequalities under nonlinear moment constraints to the stability theorem of Bianchi--Egnell \cite{BianchiEgnell1991}  through stereographic projection. Here we combine stereographic averaging with the gradient stability theorem of Figalli--Zhang \cite{FigalliZhang2022} to estimate local second moments in terms of Sobolev deficits and lower-order norms. These estimates allow us to absorb the localization errors at the optimal improved coefficient. For compactly supported functions, we obtain the stronger estimate needed when $2<p<n<3p-2$. To prove the failure of \eqref{eq:main-p} for $p>2$ and $n\geq3p-2$, we use antipodal copies of Druet's test functions. At $n=3p-2$, the curvature term is logarithmic. For $p=2$, the local moment argument may also be adapted to establish Sobolev inequalities with sharp gradient coefficients under the second- and third-order moment constraints considered in \cites{HangWang2022,Putterman2020}.

This paper is organized as follows. In Section~\ref{sec:estimates}, we establish the stereographic deficit estimates and local moment bounds. In Section~\ref{sec:endpoint}, we prove Theorem~\ref{thm:main}, including its sharpness and counterexamples for $p>2$ at and above the dimension threshold.

Throughout the paper, we write $d(x,z)$ for the spherical distance and $B_r(z)$ for a geodesic ball; balls in $\Rn$ are Euclidean. Gradients and $L^r$-norms are taken on the indicated domains. Convergence of measures is understood against continuous test functions. Constants $C>0$ may change from line to line and, unless stated otherwise, depend on the fixed parameters and localization regions. All subsequences are relabeled.

\section{Stereographic stability and local moments}\label{sec:estimates}

In this section, we establish the second moment estimates for localized
critical densities needed in the proof of Theorem~\ref{thm:main}.
We begin by averaging the Euclidean Sobolev deficit over stereographic poles.

For a pole $z\in\Sn$, set $$\Omega_z(x)=1-x\cdot z$$ and, for
$x\in\Sn\setminus\{z\}$, define
\[
b_z(x)=\frac{n-p}{p}\nabla\log\Omega_z(x)
 =-\frac{n-p}{p}\frac{z-(x\cdot z)x}{1-x\cdot z}.
 \]
In stereographic coordinates $y=\pi_z(x)$ with pole $z$, the round
metric is $\Omega_z^2\,dy^2$, where $\Omega_z=2/(1+|y|^2)$. Define
\[
 T_zf(y)=\Omega_z(\pi_z^{-1}(y))^{\frac{n-p}{p}}f(\pi_z^{-1}(y)).
\]
For $f\in W^{1,p}(\Sn)$, change of variables and the Sobolev chain
rule on compact subsets of $\Rn$ give
\begin{equation}\label{eq:stereographic}
 \begin{aligned}
 \|T_zf\|_{L^{np/(n-p)}(\Rn)}&=\|f\|_{L^{np/(n-p)}(\Sn)},\\
 \int_{\Rn}|\nabla_y T_zf|^p\,dy
 &=\int_{\Sn}|\nabla_{\Sn} f+b_zf|^p\ds.
 \end{aligned}
\end{equation}
Rotation invariance and $|b_z(x)|=(n-p)\cot(d(x,z)/2)/p$ imply
\begin{equation}\label{eq:pole-average}
 \frac1{|\Sn|}\int_{\Sn}b_z(x)\ds(z)=0,
 \qquad
 \frac1{|\Sn|}\int_{\Sn}|b_z(x)|^a\ds(z)<\infty
 \quad(0<a<n).
\end{equation}
The latter integral is independent of $x$.

We first verify that the Euclidean Sobolev inequality applies for almost
every pole. By Fubini's theorem and \eqref{eq:pole-average},
\[
 \int_{\Sn}\|b_zf\|_{L^p(\Sn)}^p\ds(z)=C(n,p)\|f\|_{L^p(\Sn)}^p<\infty.
\]
Thus $b_zf\in L^p(\Sn)$ for almost every $z$. For each such pole,
\eqref{eq:stereographic} gives $T_zf\in L^{np/(n-p)}(\Rn)$ and
$\nabla T_zf\in L^p(\Rn)$. Let $\chi_R$ be a smooth cutoff with $\chi_R=1$ on $B_R(0)$,
$\supp\chi_R\subset B_{2R}(0)$, and $|\nabla\chi_R|\leq C/R$. H\"older's inequality gives
\[
 \begin{aligned}
 &\|\nabla(\chi_RT_zf)-\nabla T_zf\|_{L^p(\Rn)}\leq\|\nabla T_zf\|_{L^p(\Rn\setminus B_R(0))}
       +C\|T_zf\|_{L^{np/(n-p)}(B_{2R}(0)\setminus B_R(0))}\to0
 \end{aligned}
\]
as $R\to\infty$. Mollifying $\chi_RT_zf$, we obtain $T_zf\in\mathring W^{1,p}(\Rn)$,
where $\mathring W^{1,p}(\Rn)$ denotes the completion of $C_c^\infty(\Rn)$
under the norm $\|\nabla\cdot\|_{L^p(\Rn)}$. The Euclidean Sobolev
inequality therefore gives
\[
 \delta_z(f):=S_{n,p}^p\|\nabla T_zf\|_{L^p(\Rn)}^p
                   -\|f\|_{L^{np/(n-p)}(\Sn)}^p\geq0
\]
for almost every $z$.

For vectors $a,b$, Taylor's inequality gives
\begin{equation}\label{eq:taylor}
 |a+b|^p\leq |a|^p+p|a|^{p-2}a\cdot b+
 \begin{cases}
 C|b|^p,&1<p\leq2,\\
 C\bigl(|a|^{p-2}|b|^2+|b|^p\bigr),&p>2.
 \end{cases}
\end{equation}
Averaging \eqref{eq:stereographic} and using
\eqref{eq:pole-average} to cancel the linear term in
\eqref{eq:taylor}, we obtain
\begin{equation}\label{eq:averaged-deficit}
 \begin{aligned}
 \frac1{|\Sn|}\int_{\Sn}\delta_z(f)\ds(z)
 &\leq S_{n,p}^p\|\nabla f\|_{L^p(\Sn)}^p-\|f\|_{L^{np/(n-p)}(\Sn)}^p\\
 &\quad+
 \begin{cases}
 C\|f\|_{L^p(\Sn)}^p,&p\leq2,\\
 C\bigl(\int_{\Sn}|\nabla f|^{p-2}|f|^2\ds+\|f\|_{L^p(\Sn)}^p\bigr),&p>2.
 \end{cases}
 \end{aligned}
\end{equation}
For $p>2$ and $\|\nabla f\|_{L^p(\Sn)}+\|f\|_{L^p(\Sn)}\leq M$, H\"older's inequality
gives
\[
 \begin{aligned}
 \int_{\Sn}|\nabla f|^{p-2}|f|^2\ds+\|f\|_{L^p(\Sn)}^p
 &\leq\|\nabla f\|_{L^p(\Sn)}^{p-2}\|f\|_{L^p(\Sn)}^2
       +\|f\|_{L^p(\Sn)}^p\leq2M^{p-2}\|f\|_{L^p(\Sn)}^2.
 \end{aligned}
\]
We may therefore choose $C_0=C(n,p)$ for $1<p\leq2$ and
$C_0=C(n,p)M^{p-2}$ for $p>2$ such that
\begin{equation}\label{eq:deficit}
 \begin{aligned}
 D(f)&:=S_{n,p}^p\|\nabla f\|_{L^p(\Sn)}^p-\|f\|_{L^{np/(n-p)}(\Sn)}^p
             +C_0\|f\|_{L^p(\Sn)}^{\min\{p,2\}},\\
 0\leq\frac1{|\Sn|}\int_{\Sn}\delta_z(f)\ds(z)&\leq D(f).
 \end{aligned}
\end{equation}
In particular, when $1<p\leq2$, \eqref{eq:deficit} gives,
for every $f\in W^{1,p}(\Sn)$,
\begin{equation}\label{eq:unconstrained}
 \|f\|_{L^{np/(n-p)}(\Sn)}^p
 \leq S_{n,p}^p\|\nabla f\|_{L^p(\Sn)}^p+C_0\|f\|_{L^p(\Sn)}^p.
\end{equation}

We next estimate the second moment of a concentrating sequence.

\begin{lemma}\label{lem:local-moment}
Let $n\geq2$ and $2n/(n+2)<p<n$. Suppose that $f_k\in W^{1,p}(\Sn)$, $f_k\geq0$, and
$\sup_k\|\nabla f_k\|_{L^p(\Sn)}<\infty$, with $\supp f_k$ contained in
a fixed cap of radius less than $\pi/4$. Let $a>0$ and let $\xi$ lie in
the interior of this cap. Assume
\begin{equation}\label{eq:local-hypotheses}
 \begin{gathered}
 f_k^{np/(n-p)}d\sigma\rightharpoonup a\delta_\xi,\qquad
 \|f_k\|_{L^p(\Sn)}\to0,\qquad d_k:=D(f_k)\to0,
 \end{gathered}
\end{equation}
where the measures converge weakly and $D$ is defined by
\eqref{eq:deficit} for this bounded family. Then there exist
$\xi_k\to\xi$ such that
\begin{equation}\label{eq:local-moment}
\frac12  \int_{\Sn}|x-\xi_k|^2f_k^{\pc}\ds= \int_{\Sn}(1-x\cdot\xi_k)f_k^{\pc}\ds
 \leq C\|f_k\|_{L^p(\Sn)}^{\min\{p,2\}}+o(d_k).
\end{equation}
\end{lemma}

\begin{proof}
Set $\gamma=\max\{2,p\}$.
Fix the normalized Talenti profile
\begin{equation}\label{eq:talenti}
 V(y)=c_{n,p}(1+|y|^{p/(p-1)})^{-\frac{n-p}{p}},
 \qquad \|V\|_{L^{np/(n-p)}(\Rn)}=1.
\end{equation}
Let $K$ be the cap containing $\supp f_k$, and choose a larger cap
$K'$ with $\overline K\subset K'$ and $|\Sn\setminus K'|>0$.
The nonnegativity of $\delta_z(f_k)$ and \eqref{eq:deficit} give
\[
 \int_{\Sn\setminus K'}\delta_z(f_k)\ds(z)
 \leq\int_{\Sn}\delta_z(f_k)\ds(z)\leq|\Sn|d_k.
\]
Hence there are $z_k\in\Sn\setminus K'$ such that
\[
 \delta_{z_k}(f_k)
 \leq\frac1{|\Sn\setminus K'|}
      \int_{\Sn\setminus K'}\delta_z(f_k)\ds(z)
 \leq\frac{|\Sn|}{|\Sn\setminus K'|}d_k=C d_k.
\]
Here $C$ is independent of $k$. For large $k$,
$T_{z_k}f_k\not\equiv0$ and $T_{z_k}f_k$ has compact support.
Hence it is not a Euclidean extremal, so $\delta_{z_k}(f_k)>0$.
By \eqref{eq:stereographic} and the definition of $\delta_{z_k}(f_k)$,
\[
 \begin{aligned}
 \|T_{z_k}f_k\|_{L^{np/(n-p)}(\Rn)}
 &=\|f_k\|_{L^{np/(n-p)}(\Sn)}\to a^{(n-p)/(np)}>0,\\
 S_{n,p}^p\|\nabla T_{z_k}f_k\|_{L^p(\Rn)}^p
 &=\|T_{z_k}f_k\|_{L^{np/(n-p)}(\Rn)}^p+\delta_{z_k}(f_k)\leq C.
 \end{aligned}
\]
Hence
\[
 \begin{aligned}
 \varepsilon_k
 &:=\frac{S_{n,p}\|\nabla T_{z_k}f_k\|_{L^p(\Rn)}}
          {\|T_{z_k}f_k\|_{L^{np/(n-p)}(\Rn)}}-1\\
 &=\left(1+\frac{\delta_{z_k}(f_k)}
                    {\|T_{z_k}f_k\|_{L^{np/(n-p)}(\Rn)}^p}\right)^{1/p}-1\\
 &\leq\frac{\delta_{z_k}(f_k)}
                {p\|T_{z_k}f_k\|_{L^{np/(n-p)}(\Rn)}^p}\\
                &\leq C d_k.
 \end{aligned}
\]
By Figalli--Zhang \cite[Theorem~1.1]{FigalliZhang2022}, there are
extremals
\[
 U_k(y)=c_k\lambda_k^{-\frac{n-p}{p}}
                   V\left(\frac{y-y_k}{\lambda_k}\right),
 \qquad c_k\neq0,\quad\lambda_k>0,\quad y_k\in\Rn,
\]
such that
\[
 \left(\frac{\|\nabla(T_{z_k}f_k-U_k)\|_{L^p(\Rn)}}{\|\nabla T_{z_k}f_k\|_{L^p(\Rn)}}\right)^\gamma
 \leq C\varepsilon_k\leq C d_k.
\]
For $W_k=T_{z_k}^{-1}U_k$, \eqref{eq:stereographic} and the Euclidean
Sobolev inequality give
\begin{equation}\label{eq:stability}
 \begin{aligned}
 \|f_k-W_k\|_{L^{np/(n-p)}(\Sn)}^\gamma
 &=\|T_{z_k}f_k-U_k\|_{L^{np/(n-p)}(\Rn)}^\gamma\\
 &\leq S_{n,p}^\gamma\|\nabla(T_{z_k}f_k-U_k)\|_{L^p(\Rn)}^\gamma\\
 &\leq C d_k.
 \end{aligned}
\end{equation}
Since $T_{z_k}f_k\geq0$ and
$\|T_{z_k}f_k\|_{L^{np/(n-p)}(\Rn)}\to a^{(n-p)/(np)}>0$,
$c_k>0$ for large $k$ and
$c_k=\|U_k\|_{L^{np/(n-p)}(\Rn)}\to a^{(n-p)/(np)}$.

Passing to a subsequence and choosing convergent coordinate frames,
assume $z_k\to z_\infty$ and
$\pi_{z_k}\to\pi_{z_\infty}$ in $C^1(\overline K)$. Set
$y_\infty=\pi_{z_\infty}(\xi)$. By \eqref{eq:stability},
$\|(T_{z_k}f_k)^{\frac{np}{n-p}}-U_k^{\frac{np}{n-p}}\|_{L^1(\Rn)}=o(1)$, and hence
$U_k^{\frac{np}{n-p}}\,dy\rightharpoonup a\delta_{y_\infty}$. For each fixed $r>0$,
\[
 \begin{aligned}
 a+o(1)=\int_{B_r(y_\infty)}U_k^{\frac{np}{n-p}}\,dy
 &\leq C\left(\frac r{\lambda_k}\right)^n,\\
 \int_{\Rn\setminus B_r(y_k)}U_k^{\frac{np}{n-p}}\,dy
 &=c_k^{\frac{np}{n-p}}\int_{|y|>r/\lambda_k}V(y)^{\frac{np}{n-p}}\,dy\to 0.
 \end{aligned}
\]
The first line gives $\limsup_k\lambda_k\leq Cr$ for every $r>0$,
so $\lambda_k\to0$; the second then gives $y_k\to y_\infty$.
Let $F_k=\pi_{z_k}^{-1}$, $\xi_k=F_k(y_k)\to\xi$, and
$\rho_k(x)^2=1-x\cdot\xi_k$. We have
$\rho_k(F_k(y))\leq C|y-y_k|$ for $y\in\Rn$ and
$C^{-1}\leq\Omega_{z_k}\leq C$ on $\overline K$.
Scaling and H\"older's inequality yield
\[
 \begin{aligned}
 c\lambda_k
 &\leq\|W_k\|_{L^p(B_{\lambda_k}(\xi_k))}\\
 &\leq\|f_k\|_{L^p(\Sn)}
       +C\lambda_k\|f_k-W_k\|_{L^{np/(n-p)}(\Sn)}\\
 &=\|f_k\|_{L^p(\Sn)}+o(\lambda_k).
 \end{aligned}
\]
Thus
\begin{equation}\label{eq:scale}
 \lambda_k^p\leq C\|f_k\|_{L^p(\Sn)}^p.
\end{equation}

If $p\leq2$, then $n>2(p-1)$ and
\[
 \begin{aligned}
 \int_{\Sn}\rho_k^2W_k^{\frac{np}{n-p}}\ds
 &=\int_{\Rn}\rho_k(F_k(y))^2U_k(y)^{\frac{np}{n-p}}\,dy\\
 &\leq C\lambda_k^2\int_{\Rn}|y|^2V(y)^{\frac{np}{n-p}}\,dy\\
& \leq C\lambda_k^2.
 \end{aligned}
\]
Using $\rho_k^2\leq2$, $np/(n-p)>2$, and \eqref{eq:stability}--\eqref{eq:scale},
we obtain
\[
 \begin{aligned}
 \int_{\Sn}\rho_k^2f_k^{\frac{np}{n-p}}\ds
 &\leq C\lambda_k^2+C\|f_k-W_k\|_{L^{np/(n-p)}(\Sn)}^{\frac{np}{n-p}}\\
 &\leq C\|f_k\|_{L^p(\Sn)}^p+C d_k^{\frac{np}{2(n-p)}}\\
& =C\|f_k\|_{L^p(\Sn)}^p+o(d_k).
 \end{aligned}
\]

If $p>2$, set $s=p/(p-2)$. Since $\frac{n-p}{p}(\frac{np}{n-p}-2)s=n+2s$, rescaling gives
\begin{equation}\label{eq:weighted-bubble}
 \begin{aligned}
 \|\rho_k^2W_k^{\frac{np}{n-p}-2}\|_{L^s(K)}
 &\leq C\bigl\||y-y_k|^2U_k^{\frac{np}{n-p}-2}\bigr\|_{L^s(\Rn)}\\
 &\leq C\left(\int_{\Rn}|y|^{2s}
 (1+|y|^{p/(p-1)})^{-n-2s}\,dy\right)^{1/s}\\
 &\leq C.
 \end{aligned}
\end{equation}
The radial integrand is $O(r^{n+2s-1})$ at zero and
$O(r^{-1-(n+2s)/(p-1)})$ at infinity.
Using $f_k^{\frac{np}{n-p}}\leq C f_k^2W_k^{\frac{np}{n-p}-2}+C|f_k-W_k|^{\frac{np}{n-p}}$ and H\"older's inequality,
we conclude from \eqref{eq:weighted-bubble} and \eqref{eq:stability} that
\[
 \begin{aligned}
 \int_{\Sn}\rho_k^2f_k^{\frac{np}{n-p}}\ds
 &\leq C\|f_k\|_{L^p(\Sn)}^2\|\rho_k^2W_k^{\frac{np}{n-p}-2}\|_{L^s(K)}+C\|f_k-W_k\|_{L^{np/(n-p)}(\Sn)}^{\frac{np}{n-p}}\\
 &\leq C\|f_k\|_{L^p(\Sn)}^2+C d_k^{\frac{n}{n-p}}\\
 &=C\|f_k\|_{L^p(\Sn)}^2+o(d_k).
 \end{aligned}
\]
\end{proof}

For $2<p<n<3p-2$, compact support yields a second moment estimate
of smaller order than the local Sobolev deficit. We also estimate the
gradient away from the concentration point to control the error
introduced by the partition of unity.

\begin{lemma}\label{lem:compact-support}
Let $2<p<n<3p-2$ and $a>0$. There exists $R_0=R_0(n,p,a)>0$ such that
the following holds for each fixed $0<R<R_0$. Suppose that $f_k\in W^{1,p}(\Sn)$, $f_k\geq0$,
$\supp f_k\subset B_R(\xi)$, and
\[
 \begin{gathered}
 f_k^{np/(n-p)}d\sigma\rightharpoonup a\delta_\xi,\quad
 \|f_k\|_{L^p(\Sn)}\to 0, \quad
 \mathcal D_k:=S_{n,p}^p\|\nabla f_k\|_{L^p(\Sn)}^p
                  -\|f_k\|_{L^{np/(n-p)}(\Sn)}^p\to 0.
 \end{gathered}
\]
Then $\mathcal D_k>0$ for large $k$, and there exist $\xi_k\to\xi$ such that,
for every fixed $\rho>0$,
\begin{equation}\label{eq:compact-moments}
 \begin{aligned}
 \int_{\Sn}(1-x\cdot\xi_k)f_k^{\pc}\ds&=o(\mathcal D_k),\\
 \int_{\Sn\setminus B_\rho(\xi)}|\nabla f_k|^p\ds
       &\leq C_{\rho,R}\mathcal D_k.
 \end{aligned}
\end{equation}
\end{lemma}

\begin{proof}
Take $R_0<\pi/8$, to be decreased below. We first derive estimates
uniform in $0<R<R_0$ and then choose $R_0$ small enough to absorb the
error terms. Set 
$$J_k=\int_{\Sn}|\nabla f_k|^{p-2}f_k^2\ds.$$ By
\eqref{eq:averaged-deficit},
\[
 0\leq\frac1{|\Sn|}\int_{\Sn}\delta_z(f_k)\ds(z)
 \leq\mathcal D_k+CJ_k+C\|f_k\|_{L^p(\Sn)}^p.
\]
By the assumptions, $\sup_k\|\nabla f_k\|_{L^p(\Sn)}<\infty$, and hence, H\"older's
inequality gives $J_k\to0$.
Since $|\{z\in\Sn:d(z,\xi)\geq\pi/2\}|=|\Sn|/2$ and
$\delta_z(f_k)\geq0$, we may choose $z_k$ with $d(z_k,\xi)\geq\pi/2$
such that
\begin{equation}\label{eq:compact-selected}
 0<\delta_{z_k}(f_k)
 \leq2\bigl(\mathcal D_k+CJ_k+C\|f_k\|_{L^p(\Sn)}^p\bigr)
 \rightarrow0.
\end{equation}
As before, strict positivity follows because a nonzero compactly supported
function is not a Euclidean extremal.

Set $g_k=T_{z_k}f_k$. The selected charts and their conformal factors are uniformly
bounded on $\supp f_k$, independently of small $R$.
By the gradient stability theorem of Figalli--Zhang
\cite[Theorem~1.1]{FigalliZhang2022}, choose extremals
$U_k$ of the form used in the proof of Lemma~\ref{lem:local-moment},
with profile \eqref{eq:talenti}, centers $y_k$ and scales
$\lambda_k\to0$, such that % $e_k=g_k-U_k$ satisfies
\begin{equation}\label{eq:compact-stability}
 \|\nabla (g_k-U_k)\|_{L^p(\Rn)}^p
       +\|g_k-U_k\|_{L^{np/(n-p)}(\Rn)}^p\leq C\delta_{z_k}(f_k).
\end{equation}
Since $g_k\geq0$ and
$\|g_k\|_{L^{np/(n-p)}(\Rn)}\to a^{(n-p)/(np)}>0$,
\eqref{eq:compact-stability} implies $c_k>0$ for large $k$.
By \eqref{eq:talenti} and \eqref{eq:stereographic},
\[
 \left|c_k-\|f_k\|_{L^{np/(n-p)}(\Sn)}\right|
 \leq\|g_k-U_k\|_{L^{np/(n-p)}(\Rn)}\to0.
\]
Thus $c_k\to a^{(n-p)/(np)}>0$. Moreover,
\eqref{eq:compact-stability} gives
$\|g_k^{np/(n-p)}-U_k^{np/(n-p)}\|_{L^1(\Rn)}=o(1)$.
For each fixed $r>0$, the uniform Lipschitz bounds for the charts,
the concentration assumption, and $\lambda_k\to0$ yield
\[
 \begin{aligned}
 \int_{\Rn\setminus B_r(\pi_{z_k}(\xi))}
       U_k^{\frac{np}{n-p}}\,dy
 &\leq\int_{\Sn\setminus B_{r/C}(\xi)}
       f_k^{\frac{np}{n-p}}\ds+o(1)=o(1),\\
 \int_{B_r(y_k)}U_k^{\frac{np}{n-p}}\,dy
 &=c_k^{\frac{np}{n-p}}\int_{|y|<r/\lambda_k}
       V(y)^{\frac{np}{n-p}}\,dy=a+o(1).
 \end{aligned}
\]
If $|y_k-\pi_{z_k}(\xi)|\geq2r$ along a subsequence, then
$B_r(y_k)\subset\Rn\setminus B_r(\pi_{z_k}(\xi))$,
contradicting $a>0$. Hence $|y_k-\pi_{z_k}(\xi)|\to0$. For large $k$,
$\supp g_k\subset B_{CR}(y_k)$, with $C$ independent of $R$.
Then the explicit
decay of $\nabla U_k$ gives
\begin{equation}\label{eq:compact-tail}
 c\lambda_k^\alpha R^{-\alpha}
 \leq\int_{|y-y_k|>CR}|\nabla U_k|^p\,dy
 \leq\|\nabla (g_k-U_k )\|_{L^p(\Rn)}^p\leq C\delta_{z_k}(f_k),
\end{equation}
where $\alpha=(n-p)/(p-1)<2$. 
By \eqref{eq:compact-stability} and
\eqref{eq:compact-tail}, radial integration and H\"older's inequality give
\begin{equation}\label{eq:compact-profile}
 \begin{aligned}
 \int_{B_{CR}(y_k)}U_k^p\,dy
     &\leq C\lambda_k^\alpha R^{p-\alpha}
      \leq CR^p\delta_{z_k}(f_k),\\
 \int_{B_{CR}(y_k)}|\nabla U_k|^{p-2}U_k^2\,dy
     &\leq C\lambda_k^\alpha R^{2-\alpha}
      \leq CR^2\delta_{z_k}(f_k),\\
 \int_{B_{CR}(y_k)}|g_k-U_k|^p\,dy
     &\leq CR^p\delta_{z_k}(f_k).
 \end{aligned}
\end{equation}

Since $\supp g_k\subset B_{CR}(y_k)$, \eqref{eq:compact-profile} gives
\begin{equation}\label{eq:compact-Lp}
 \|g_k\|_{L^p(\Rn)}^p\leq C\int_{B_{CR}(y_k)}\bigl(U_k^p+|g_k-U_k|^p\bigr)\,dy\leq CR^p\delta_{z_k}(f_k).
\end{equation}
Hence H\"older's inequality and \eqref{eq:compact-stability} imply
\[
 \begin{aligned}
 \int_{\Rn}|\nabla(g_k-U_k)|^{p-2}g_k^2\,dy
 &\leq\|\nabla(g_k-U_k)\|_{L^p(\Rn)}^{p-2}
       \|g_k\|_{L^p(\Rn)}^2\\
 &\leq C\delta_{z_k}(f_k)^{(p-2)/p}
       \bigl(R^p\delta_{z_k}(f_k)\bigr)^{2/p}\\
 &=CR^2\delta_{z_k}(f_k).
 \end{aligned}
\]
To estimate $\displaystyle\int_{\Rn}|\nabla U_k|^{p-2}g_k^2\,dy$, split $\supp g_k$
into $\{g_k\leq2U_k\}$ and $\{g_k>2U_k\}$. On the first set,
$|\nabla U_k|^{p-2}g_k^2\leq4|\nabla U_k|^{p-2}U_k^2$. On the second set, $g_k-U_k>U_k$ and
$g_k\leq2(g_k-U_k)$. Since $|\nabla U_k|/U_k\leq\alpha/|y-y_k|$,
\[
 \begin{aligned}
 |\nabla U_k|^{p-2}g_k^2
 &\leq C|y-y_k|^{2-p}U_k^{p-2}|g_k-U_k|^2\leq C|y-y_k|^{2-p}|g_k-U_k|^p.
 \end{aligned}
\]
Since $1<p<n$ and $g_k-U_k\in\mathring W^{1,p}(\Rn)$, Hardy's
inequality and \eqref{eq:compact-stability} give
\[
 \begin{aligned}
 \int_{B_{CR}(y_k)}|y-y_k|^{2-p}|g_k-U_k|^p\,dy
 &\leq CR^2\int_{\Rn}\frac{|g_k-U_k|^p}{|y-y_k|^p}\,dy\\
 &\leq CR^2\|\nabla(g_k-U_k)\|_{L^p(\Rn)}^p\\
 &\leq CR^2\delta_{z_k}(f_k).
 \end{aligned}
\]
Together with \eqref{eq:compact-profile}, these estimates yield
$$\int_{\Rn}|\nabla U_k|^{p-2}g_k^2\,dy\leq CR^2\delta_{z_k}(f_k).$$
Using $|\nabla g_k|^{p-2}\leq C\bigl(|\nabla U_k|^{p-2}
+|\nabla(g_k-U_k)|^{p-2}\bigr)$, we conclude that
\begin{equation}\label{eq:compact-gradient}
 \begin{aligned}
% \|g_k\|_{L^p(\Rn)}^p&\leq CR^p\delta_{z_k}(f_k),\\
 \int_{\Rn}|\nabla g_k|^{p-2}g_k^2\,dy&\leq CR^2\delta_{z_k}(f_k).
 \end{aligned}
\end{equation}

Writing $\Omega=\Omega_{z_k}$ in the selected chart, we have
\[
 \begin{aligned}
 \|f_k\|_{L^p(\Sn)}^p&=\int_{\Rn}\Omega^p g_k^p\,dy\leq CR^p\delta_{z_k}(f_k),\\
 J_k&=\int_{\Rn}\Omega^2
   \left|\nabla g_k-\frac{n-p}{p}g_k\nabla\log\Omega\right|^{p-2}
   g_k^2\,dy\\
 &\leq C\int_{\mathbb R^n}
 \left(|\nabla_y g_k|^{p-2}g_k^2+g_k^p\right)\,dy\\
& \leq C(R^2+R^p)\delta_{z_k}(f_k).
 \end{aligned}
\]
Hence, \eqref{eq:compact-selected}
gives $\delta_{z_k}(f_k)\leq2\mathcal D_k+CR^2\delta_{z_k}(f_k)$. Choose $R_0$ so that
$CR_0^2\leq1/2$. For each fixed $0<R<R_0$, absorption yields
$\delta_{z_k}(f_k)\leq4\mathcal D_k$, hence $$\mathcal D_k>0$$ for large $k$.

For fixed $R$, \eqref{eq:compact-tail} gives
$\lambda_k^\alpha\leq CR^\alpha\delta_{z_k}(f_k)$. Since
$2-\alpha>0$ and $n/(p-1)-\alpha=p/(p-1)>0$, a direct calculation yields
\[
 \begin{aligned}
 \int_{B_{CR}(y_k)}|y-y_k|^2U_k^{\pc}\,dy
 \leq C_R\bigl(\lambda_k^2+\lambda_k^{n/(p-1)}\bigr)
                         (1+|\log\lambda_k|)=o(\lambda_k^\alpha)=o(\delta_{z_k}(f_k)).
 \end{aligned}
\]
Set $\xi_k=\pi_{z_k}^{-1}(y_k)$. Since $\Omega\leq2$, the inverse
charts are $2$-Lipschitz, and
$d(\xi_k,\xi)\leq2|y_k-\pi_{z_k}(\xi)|\to0$.
Using $1-x\cdot\xi_k=\frac12|x-\xi_k|^2$,
\eqref{eq:stereographic}, and \eqref{eq:compact-stability}, we obtain
\[
 \begin{aligned}
 \int_{\Sn}(1-x\cdot\xi_k)f_k^{\pc}\ds
 &=\int_{B_{CR}(y_k)}
       \bigl(1-\pi_{z_k}^{-1}(y)\cdot\xi_k\bigr)g_k^{\pc}\,dy\\
 &\leq C\int_{B_{CR}(y_k)}|y-y_k|^2U_k^{\pc}\,dy
       +CR^2\|g_k-U_k\|_{L^{np/(n-p)}(\Rn)}^{np/(n-p)}\\
 &\leq o(\delta_{z_k}(f_k))
       +CR^2\delta_{z_k}(f_k)^{n/(n-p)}\\
&=o(\mathcal D_k),
 \end{aligned}
\]
where $n/(n-p)>1$ and $\delta_{z_k}(f_k)\leq4\mathcal D_k$.

For fixed $\rho>0$ and large $k$, $d(\xi_k,\xi)<\rho/2$ gives
\[
 \pi_{z_k}(\supp f_k\setminus B_\rho(\xi))
 \subset\{y\in\Rn:|y-y_k|\geq\rho/4\}.
\]
By \eqref{eq:talenti}, \eqref{eq:compact-stability}, and
\eqref{eq:compact-tail},
\[
 \begin{aligned}
 \int_{|y-y_k|\geq\rho/4}|\nabla g_k|^p\,dy
 &\leq C\int_{|y-y_k|\geq\rho/4}|\nabla U_k|^p\,dy
       +C\|\nabla(g_k-U_k)\|_{L^p(\Rn)}^p\\
 &\leq C_\rho\lambda_k^\alpha+C\delta_{z_k}(f_k)\\
 &\leq C_{\rho,R}\delta_{z_k}(f_k).
 \end{aligned}
\]
Then it follows from \eqref{eq:compact-Lp} that
\[
 \begin{aligned}
 \int_{\Sn\setminus B_\rho(\xi)}|\nabla f_k|^p\ds
 &\leq\int_{|y-y_k|\geq\rho/4}
       \left|\nabla g_k-\frac{n-p}{p}g_k\nabla\log\Omega\right|^p\,dy\\
 &\leq C\int_{|y-y_k|\geq\rho/4}|\nabla g_k|^p\,dy
       +C\|g_k\|_{L^p(\Rn)}^p\\
 &\leq C_{\rho,R}\delta_{z_k}(f_k)\\
  &\leq C_{\rho,R}\mathcal D_k.
 \end{aligned}
\]
\end{proof}

\section{Proof of the main theorem}\label{sec:endpoint}

\begin{proof}[Proof of Theorem~\ref{thm:main}]
Replacing $u$ by $|u|$, we may assume $u\geq0$.

\medskip
\noindent\textit{Step 1. Averaging for $1<p<2n/(n+1)$.}
We prove \eqref{eq:main-p}. Set $d\widehat\sigma(z)=d\sigma(z)/|\Sn|$ and
\[
 \psi_z(x)=\operatorname{sgn}(x\cdot z)
                 |x\cdot z|^{\frac{n-p}{np}}.
\]
By rotation invariance, the constants
\begin{equation}\label{eq:root-integrals}
 A_0=\int_{\Sn}|\psi_z(x)|^p\,d\widehat\sigma(z)>0,
 \qquad
 A_1=\int_{\Sn}|\nabla_x\psi_z(x)|^p\,d\widehat\sigma(z)<\infty
\end{equation}
are independent of $x$. Here
$|\nabla_x\psi_z|^p\leq C|x\cdot z|^{\frac{n-p}{n}-p}$, and
$\frac{n-p}{n}-p>-1$ precisely when $p<2n/(n+1)$.

For $\delta>0$ and $t=x\cdot z$, set
\[
 \psi_{z,\delta}(x)=t(t^2+\delta^2)^{\frac12(\frac{n-p}{np}-1)}.
\]
Since $|\psi_{z,\delta}|\leq|t|^{(n-p)/(np)}$ and
$|\partial_t\psi_{z,\delta}|\leq C|t|^{(n-p)/(np)-1}$ for $t\ne0$,
\eqref{eq:root-integrals} and the dominated convergence theorem give, as $\delta\downarrow0$,
\[
 \begin{aligned}
 \psi_{z,\delta}u&\longrightarrow\psi_z u,\\
 \nabla_x(\psi_{z,\delta}u)
 &\longrightarrow\psi_z\nabla u+u\nabla_x\psi_z
 \end{aligned}
 \qquad\text{in }L^p(\Sn\times\Sn).
\]
By Fubini's theorem, $\psi_z u\in W^{1,p}(\Sn)$ for almost every $z$.
The assumption \eqref{eq:constraints} leads to
 \begin{align*}
 \int_{\Sn}(\psi_z^+u)^{\pc}\ds&=\int_{\Sn}(\psi_z^-u)^{\pc}\ds,\\
 \|\psi_z^+u\|_{L^{np/(n-p)}(\Sn)}^p+\|\psi_z^-u\|_{L^{np/(n-p)}(\Sn)}^p&=2^{p/n}\|\psi_z u\|_{L^{np/(n-p)}(\Sn)}^p,\\
 |\nabla_x(\psi_z^+u)|^p+|\nabla_x(\psi_z^-u)|^p&=|\nabla_x(\psi_z u)|^p\quad\text{a.e.}
 \end{align*}

Applying \eqref{eq:unconstrained} to $\psi_z^\pm u$ and adding, we obtain
\[
 \|\psi_z u\|_{L^{np/(n-p)}(\Sn)}^p
 \leq 2^{-p/n}\bigl(S_{n,p}^p\|\nabla_x(\psi_z u)\|_{L^p(\Sn)}^p+C_0\|\psi_z u\|_{L^p(\Sn)}^p\bigr).
\]
Minkowski's inequality, with exponent $n/(n-p)>1$, yields
 \begin{align}
 A_0\|u\|_{L^{np/(n-p)}(\Sn)}^p
 &=\left\|\int_{\Sn}|\psi_z u|^p\,d\widehat\sigma(z)\right\|_{L^{n/(n-p)}(\Sn)}\nonumber\\
 &\leq\int_{\Sn}\|\psi_z u\|_{L^{np/(n-p)}(\Sn)}^p\,d\widehat\sigma(z)\nonumber\\
 &\leq 2^{-p/n}\left(S_{n,p}^p\int_{\Sn} \|\nabla_x(\psi_z u)\|_{L^p(\Sn)}^p\,d\widehat\sigma(z)+C_0\int_{\Sn} \|\psi_z u\|_{L^p(\Sn)}^p\,d\widehat\sigma(z)\right)\nonumber\\
  &\leq 2^{-p/n}\left(S_{n,p}^p\int_{\Sn} \|\nabla_x(\psi_z u)\|_{L^p(\Sn)}^p\,d\widehat\sigma(z)+C_0A_0\|u\|_{L^p(\Sn)}^p\right).\label{eq:root-inequality}
 \end{align}
Since $(n-p)/n-1>-1$,
\[
 \int_{\Sn}|\psi_z|^{p-2}\psi_z\nabla_x\psi_z\,d\widehat\sigma(z)=\frac{1}{p} \int_{\Sn}\nabla_x|\psi_z|^p\,d\widehat\sigma(z) =\frac1p\nabla_x A_0=0.
\]
Taking $a=\psi_z\nabla u$ and $b=u\nabla_x\psi_z$ in
\eqref{eq:taylor}, we obtain
\[
\int_{\mathbb S^n}|\nabla_x(\psi_z u)(x)|^p\,d\widehat\sigma(z)
\leq A_0|\nabla u(x)|^p+CA_1u(x)^p.
\]
Integrating over $x$ and using Tonelli's theorem, we have
\[
 \int_{\Sn}\|\nabla_x(\psi_z u)\|_{L^p(\Sn)}^p\,d\widehat\sigma(z)\leq A_0\|\nabla u\|_{L^p(\Sn)}^p+C A_1\|u\|_{L^p(\Sn)}^p.
\]
Then it follows from \eqref{eq:root-inequality} that
\[
 \|u\|_{L^{np/(n-p)}(\Sn)}^p\leq2^{-p/n}S_{n,p}^p\|\nabla u\|_{L^p(\Sn)}^p+2^{-p/n}\left(C_0+\frac{CS_{n,p}^pA_1}{A_0}\right)\|u\|_{L^p(\Sn)}^p,
\]
which is \eqref{eq:main-p}.

\medskip
\noindent\textit{Step 2. Two-point concentration and deficit absorption.}
Assume $2n/(n+1)\leq p<n$. If \eqref{eq:main-p} fails for $p\leq2$,
or \eqref{eq:main-square} fails for $p>2$, then there exists a sequence of nonnegative functions
$\{u_k\}$ satisfying \eqref{eq:constraints},
\begin{equation}\label{eq:failure}
 \begin{gathered}
 \|u_k\|_{L^{np/(n-p)}(\Sn)}=1,
 \quad \ell_k=\|u_k\|_{L^p(\Sn)}^{\min\{p,2\}},\quad 1>2^{-p/n}S_{n,p}^p\|\nabla u_k\|_{L^p(\Sn)}^p+k\ell_k,
 \end{gathered}
\end{equation}
where we use $t\leq t^{\frac{2}{p}}$ for $0\leq t\leq1$ and  $p>2$  in the last inequality.
Then $\sup_k\|\nabla u_k\|_{L^p(\Sn)}<\infty$,
$u_k\to0$ in $L^p(\Sn)$, and $u_k\rightharpoonup0$ in $W^{1,p}(\Sn)$.
By the concentration compactness in Lions \cite{Lions1985}, in the form stated by Hang--Wang
\cite[Theorem~2.1]{HangWang2022}, after passing to a subsequence,
\[
 |\nabla u_k|^p d\sigma\rightharpoonup\mu,
 \qquad u_k^{np/(n-p)}d\sigma\rightharpoonup\sum_j a_j\delta_{z_j},\qquad
 a_j^{\frac{n-p}{n}}\leq S_{n,p}^p\mu(\{z_j\}),
\]
weakly as measures on $\Sn$. Then
\[
 \sum_j a_j=1,\qquad \sum_j a_jz_j=0,\qquad
 a_j=\left|\sum_{i\ne j}a_i z_i\right|\leq1-a_j.
\]
Hence $a_j\leq1/2$, and \eqref{eq:failure} implies
\[
 2^{p/n}=2^{p/n}\sum_j a_j
 \leq\sum_j a_j^{\frac{n-p}{n}}
 \leq S_{n,p}^p\mu(\Sn)\leq2^{p/n}.
\]
Equality forces $a_j\in\{0,1/2\}$. Since $\sum_j a_j=1$ and
$\sum_j a_jz_j=0$, after relabeling, $a_1=a_2=1/2$ and
$z_2=-z_1$, with no other positive $a_j$. Moreover,
\begin{equation}\label{eq:energy-limit}
 S_{n,p}^p\|\nabla u_k\|_{L^p(\Sn)}^p\longrightarrow2^{p/n}.
\end{equation}

Choose $0<R<\pi/4$ and nonnegative Lipschitz functions $\chi_i$ such that
$\sum_{i=0}^2\chi_i^p=1$ and, for $i=1,2$,
$\supp\chi_i\subset B_R(z_i)$ and $\chi_i=1$ on $B_{R/2}(z_i)$.
Since $\sum_i\chi_i^{p-1}\nabla\chi_i
=\frac1p\nabla(\sum_i\chi_i^p)=0$, \eqref{eq:taylor} gives
\begin{equation}\label{eq:localization}
 \sum_{i=0}^2\|\nabla(\chi_i u_k)\|_{L^p(\Sn)}^p
 \leq\|\nabla u_k\|_{L^p(\Sn)}^p+C\ell_k,
\end{equation}
where for $p>2$ we used that
\[
 \begin{aligned}
 \int_{\Sn}|\nabla u_k|^{p-2}u_k^2\ds+\|u_k\|_{L^p(\Sn)}^p
 &\leq\|\nabla u_k\|_{L^p(\Sn)}^{p-2}\|u_k\|_{L^p(\Sn)}^2
       +\|u_k\|_{L^p(\Sn)}^p\leq C\ell_k.
 \end{aligned}
\]

Set $f_{i,k}=\chi_i u_k$ and
$m_{i,k}=\int_{\Sn}f_{i,k}^{np/(n-p)}\ds$.
Fix $C_0$ in \eqref{eq:deficit} for this bounded family, and set
\[
 \begin{aligned}
 d_{i,k}&=D(f_{i,k})\quad(i=1,2),\\
 d_{0,k}&=S_{n,p}^p\|\nabla f_{0,k}\|_{L^p(\Sn)}^p
                      +C_0\|f_{0,k}\|_{L^p(\Sn)}^{\min\{p,2\}},\\
 d_k&=\sum_{i=0}^2d_{i,k}.
 \end{aligned}
\]
Then $d_{i,k}\geq0$, $m_{0,k}^{(n-p)/n}\leq d_{0,k}$ (from \eqref{eq:deficit}), and
$\sum_{i=0}^2\|f_{i,k}\|_{L^p(\Sn)}^{\min\{p,2\}}\leq 3\ell_k$.
By \eqref{eq:localization},
\begin{equation}\label{eq:localized-energy}
 S_{n,p}^p\|\nabla u_k\|_{L^p(\Sn)}^p
 \geq\sum_{i=1}^2 m_{i,k}^{\frac{n-p}{n}}+d_k-C\ell_k,
 \quad\mbox{and}\quad d_k\rightarrow0,
\end{equation}
where we additionally used $m_{i,k}\to1/2$ for $i=1,2$,
$\ell_k\to0$, and \eqref{eq:energy-limit} for the convergence of $d_k$.

Since $p\geq2n/(n+1)>2n/(n+2)$, Lemma~\ref{lem:local-moment}
gives $\xi_{i,k}\to z_i$ such that
\begin{equation}\label{eq:sum-variance}
 \sum_{i=1}^2\int_{\Sn}(1-x\cdot\xi_{i,k})f_{i,k}^{\pc}\ds
 \leq C\ell_k+o(d_k).
\end{equation}
Let $T=\Sn\setminus(B_{R/2}(z_1)\cup B_{R/2}(z_2))$ and
$\tau_k=1-m_{1,k}-m_{2,k}\geq0$.
For $i=1,2$ and large $k$, $1-x\cdot\xi_{i,k}\geq C^{-1}$ on $T$.
Since $\sum_{i=1}^2\chi_i^{np/(n-p)}\leq1$ and
$u_k^{np/(n-p)}=(\sum_{i=0}^2 f_{i,k}^p)^{n/(n-p)}
\leq C\sum_{i=0}^2 f_{i,k}^{np/(n-p)}$, it follows from \eqref{eq:sum-variance} that
\begin{equation}\label{eq:mass-loss}
 \begin{aligned}
 \tau_k
 &=\int_T\left(1-\chi_1^{np/(n-p)} -\chi_2^{np/(n-p)}\right)u_k^{\pc}\ds\\
 &\leq C\sum_{i=1}^2\int_T f_{i,k}^{\pc}\ds+Cm_{0,k}\\
 &\leq C\ell_k+o(d_k)+Cd_{0,k}^{\frac n{n-p}}\\
 &=C\ell_k+o(d_k).
 \end{aligned}
\end{equation}

By \eqref{eq:constraints} and $0\leq f_{i,k}\leq u_k$,
\[
 \begin{aligned}
 m_{i,k}
 &=\frac12\int_{\Sn}(1+x\cdot\xi_{i,k})f_{i,k}^{\pc}\ds+\frac12\int_{\Sn}(1-x\cdot\xi_{i,k})f_{i,k}^{\pc}\ds\\
 &\leq\frac12\int_{\Sn}(1+x\cdot\xi_{i,k})u_k^{\pc}\ds+\frac12\int_{\Sn}(1-x\cdot\xi_{i,k})f_{i,k}^{\pc}\ds\\
 &\leq\frac12+C\ell_k+o(d_k).
 \end{aligned}
\]
Hence, by \eqref{eq:mass-loss},
\[
 \begin{aligned}
 \sum_{i=1}^2m_{i,k}^{\frac{n-p}{n}}
 &=\sum_{i=1}^2m_{i,k}m_{i,k}^{-p/n}\geq\left(\frac12+C\ell_k+o(d_k)\right)^{-p/n}(1-\tau_k)\geq2^{p/n}-C\ell_k-o(d_k).
 \end{aligned}
\]
Substituting into \eqref{eq:localized-energy}, we obtain
\[
 1-2^{-p/n}S_{n,p}^p\|\nabla u_k\|_{L^p(\Sn)}^p
 \leq C\ell_k-2^{-p/n}d_k+o(d_k)\leq C\ell_k
\]
for large $k$, contradicting \eqref{eq:failure}.

\medskip
\noindent\textit{Step 3. The $p$-power estimate for $2<p<n<3p-2$.}
Suppose \eqref{eq:main-p} fails, and choose $u_k$ as in
\eqref{eq:failure}, with $\ell_k=\|u_k\|_{L^p(\Sn)}^p$.
Step~2 gives $u_k^{np/(n-p)}d\sigma
\rightharpoonup\frac12(\delta_{z_1}+\delta_{z_2})$,
$z_2=-z_1$, and \eqref{eq:energy-limit}.
Choose $0<R<\min\{\pi/4,R_0(n,p,1/2)\}$ in Lemma~\ref{lem:compact-support}
and use $\chi_i$, $f_{i,k}$, $m_{i,k}$ and $T$ as in Step~2. Set
\[
 \mathcal D_{i,k}=S_{n,p}^p\|\nabla f_{i,k}\|_{L^p(\Sn)}^p
                    -m_{i,k}^{\frac{n-p}{n}}\qquad(i=1,2).
\]
By \eqref{eq:taylor}, \eqref{eq:deficit}, and \eqref{eq:energy-limit},
\[
\begin{aligned}
\mathcal D_{i,k}&\geq-C\ell_k^{2/p}\qquad(i=1,2),\\
\mathcal D_{1,k}+\mathcal D_{2,k}+S_{n,p}^p\|\nabla f_{0,k}\|_{L^p(\Sn)}^p&\leq S_{n,p}^p\|\nabla u_k\|_{L^p(\Sn)}^p-\sum_{i=1}^2m_{i,k}^{\frac{n-p}{n}}+C_R(\ell_k^{2/p}+\ell_k)\\
&=o(1).
 \end{aligned}
\]
Thus $\mathcal D_{i,k}\to0$ for $i=1,2$ and
$\|\nabla f_{0,k}\|_{L^p(\Sn)}\to0$.
Lemma~\ref{lem:compact-support} gives $\xi_{i,k}\to z_i$ and
$\mathcal D_{i,k}>0$ for large $k$. In this step, set
\[
 d_k=\mathcal D_{1,k}+\mathcal D_{2,k}
                  +S_{n,p}^p\|\nabla f_{0,k}\|_{L^p(\Sn)}^p\longrightarrow0.
\]

Since $\nabla\chi_i=0$ outside $T$ and
$\sum_i\chi_i^{p-1}\nabla\chi_i=0$, \eqref{eq:taylor} and
Young's inequality give, for every fixed $\eta>0$,
\[
 \begin{aligned}
 \sum_{i=0}^2\|\nabla f_{i,k}\|_{L^p(\Sn)}^p
 &\leq\|\nabla u_k\|_{L^p(\Sn)}^p
       +C_R\int_T|\nabla u_k|^{p-2}u_k^2\ds+C_R\ell_k\\
 &\leq\|\nabla u_k\|_{L^p(\Sn)}^p
       +\eta\int_T|\nabla u_k|^p\ds+C_{\eta,R}\ell_k.
 \end{aligned}
\]
Moreover, $\nabla u_k=\sum_{i=0}^2\chi_i^{p-1}\nabla f_{i,k}$ and
$|\nabla u_k|^p\leq\sum_{i=0}^2|\nabla f_{i,k}|^p$ by H\"older's inequality.
Hence \eqref{eq:compact-moments} gives
\[
 \begin{aligned}
 \int_T|\nabla u_k|^p\ds
 &\leq\sum_{i=1}^2\int_{\Sn\setminus B_{R/2}(z_i)}|\nabla f_{i,k}|^p\ds+\|\nabla f_{0,k}\|_{L^p(\Sn)}^p\leq C_Rd_k.
 \end{aligned}
\]
Fix $\eta>0$ so that $S_{n,p}^p\eta C_R\leq1/2$. Then
\begin{equation}\label{eq:refined-energy}
 S_{n,p}^p\|\nabla u_k\|_{L^p(\Sn)}^p
 \geq\sum_{i=1}^2m_{i,k}^{\frac{n-p}{n}}
                       +\frac12d_k-C\ell_k.
\end{equation}

By \eqref{eq:compact-moments} and the usual Sobolev inequality,
\[
 \begin{aligned}
 \sum_{i=1}^2\int_T f_{i,k}^{\pc}\ds
 &\leq C_R\sum_{i=1}^2\int_{\Sn}
           (1-x\cdot\xi_{i,k})f_{i,k}^{\pc}\ds=o(d_k),\\
 m_{0,k}^{(n-p)/n}
 &\leq C\bigl(\|\nabla f_{0,k}\|_{L^p(\Sn)}^p
                    +\|f_{0,k}\|_{L^p(\Sn)}^p\bigr)
 \leq C(d_k+\ell_k).
 \end{aligned}
\]
Consequently,
\[
 \begin{aligned}
 \tau_k:=1-m_{1,k}-m_{2,k}
 &\leq C\sum_{i=1}^2\int_Tf_{i,k}^{\pc}\ds+Cm_{0,k}\\
 &\leq o(d_k)+C(d_k+\ell_k)^{\frac n{n-p}}\\
&=o(d_k+\ell_k).
 \end{aligned}
\]
As in Step~2, \eqref{eq:constraints} gives
\[
 m_{i,k}\leq\frac12
       +\frac12\int_{\Sn}(1-x\cdot\xi_{i,k})f_{i,k}^{\pc}\ds
       =\frac12+o(d_k),\qquad i=1,2.
\]
Therefore,
\[
 \sum_{i=1}^2m_{i,k}^{\frac{n-p}{n}}
 \geq\bigl(\tfrac12+o(d_k)\bigr)^{-p/n}(1-\tau_k)
 \geq2^{p/n}-o(d_k+\ell_k).
\]
Combining this with \eqref{eq:refined-energy}, we obtain
\[
 \begin{aligned}
 S_{n,p}^p\|\nabla u_k\|_{L^p(\Sn)}^p
 &\geq2^{p/n}+\frac12d_k-C\ell_k-o(d_k+\ell_k)\geq2^{p/n}-C\ell_k,
 \end{aligned}
\]
contradicting \eqref{eq:failure}.

\medskip
\noindent\textit{Step 4. Sharpness of the gradient coefficients.}
Fix $1<p<n$, $z\in\Sn$ and $0<R<\pi/4$. Choose a smooth radial
cutoff $\eta\geq0$ with $\supp\eta\subset[0,R)$ and $\eta=1$ on $[0,R/2]$.
For $V$ in \eqref{eq:talenti}, set
\[
 w_\lambda(r)=c_\lambda\eta(r)\lambda^{-\frac{n-p}{p}}V(r/\lambda),
 \qquad r=d(x,z),
\]
where $c_\lambda>0$ satisfies
$\int_{\Sn}w_\lambda^{np/(n-p)}\ds=1$.

As $\lambda\downarrow0$, scaling gives $c_\lambda\to1$ and
$w_\lambda^{np/(n-p)}d\sigma\rightharpoonup\delta_z$.
In geodesic polar coordinates,
\[
 \begin{aligned}
 &\|\nabla w_\lambda\|_{L^p(\Sn)}^p=c_\lambda^p\int_{|y|<R/(2\lambda)}|\nabla V(y)|^p
       \left(\frac{\sin(\lambda|y|)}{\lambda|y|}\right)^{n-1}\,dy
       +O\bigl(\lambda^{(n-p)/(p-1)}\bigr)\rightarrow S_{n,p}^{-p}.
 \end{aligned}
\]
Thus $w_\lambda\rightharpoonup0$ in $W^{1,p}(\Sn)$, and compactness of
$W^{1,p}(\Sn)\hookrightarrow L^p(\Sn)$ gives
$\|w_\lambda\|_{L^p(\Sn)}\to0$. Set
\[
 u_\lambda(x)=2^{-\frac{n-p}{np}}
    \bigl(w_\lambda(d(x,z))+w_\lambda(d(x,-z))\bigr).
\]
Since $u_\lambda(-x)=u_\lambda(x)$,
$\int_{\Sn}x_j u_\lambda^{\pc}\ds=0$ for $j=1,\ldots,n+1$.
The supports in $B_R(z)$ and $B_R(-z)$ are disjoint, so
 \begin{align*}
 &\|u_\lambda\|_{L^{np/(n-p)}(\Sn)}=1,\qquad
 \|u_\lambda\|_{L^p(\Sn)}=2^{1/n}\|w_\lambda\|_{L^p(\Sn)}\rightarrow0,\\
& \|\nabla u_\lambda\|_{L^p(\Sn)}^p=2^{p/n}\|\nabla w_\lambda\|_{L^p(\Sn)}^p\rightarrow2^{p/n}S_{n,p}^{-p}.
 \end{align*}
Consequently,
\[
 \begin{aligned}
 \frac{\|u_\lambda\|_{L^{np/(n-p)}(\Sn)}^p}{\|\nabla u_\lambda\|_{L^p(\Sn)}^p}
 \rightarrow2^{-p/n}S_{n,p}^p,
\qquad \frac{\|u_\lambda\|_{L^{np/(n-p)}(\Sn)}^2}{\|\nabla u_\lambda\|_{L^p(\Sn)}^2}
\rightarrow2^{-2/n}S_{n,p}^2,
 \end{aligned}
\]
which proves sharpness in \eqref{eq:main-p} and \eqref{eq:main-square}.

\medskip
\noindent\textit{Step 5. Failure of the $p$-power inequality for $p>2$ and $n\geq3p-2$.}
With $V$ and $\eta$ as in Step~4, set
\[
 v_\lambda(r)=\eta(r)\lambda^{-\frac{n-p}{p}}V(r/\lambda).
\]
For $v_\lambda(d(x,z))$, write
\[
 \begin{aligned}
 N_\lambda=\int_{\Sn}v_\lambda^{\pc}\ds,\qquad
 E_\lambda=\int_{\Sn}|\nabla v_\lambda|^p\ds,\qquad
 P_\lambda=\int_{\Sn}v_\lambda^p\ds.
 \end{aligned}
\]
The functions $u_\lambda$ from Step~4 satisfy \eqref{eq:constraints}, and
\begin{equation}\label{eq:counterexample-ratio}
 \frac{\|u_\lambda\|_{L^{np/(n-p)}(\Sn)}^p
       -2^{-p/n}S_{n,p}^p\|\nabla u_\lambda\|_{L^p(\Sn)}^p}
      {\|u_\lambda\|_{L^p(\Sn)}^p}
 =\frac{N_\lambda^{(n-p)/n}-S_{n,p}^pE_\lambda}
      {2^{p/n}P_\lambda}.
\end{equation}

\emph{The case $n>3p-2$.} Let $\alpha=(n-p)/(p-1)>2$. Then
\[
 L_2=\int_{\Rn}|y|^2V^{\pc}\,dy<\infty,\qquad
 G_2=\int_{\Rn}|y|^2|\nabla V|^p\,dy<\infty.
\]
Since $|\nabla v_\lambda|=|v_\lambda'(r)|$ and
$(\sin r/r)^{n-1}=1-(n-1)r^2/6+O(r^4)$, the dominated convergence theorem and the
cutoff error $O(\lambda^\alpha)=o(\lambda^2)$ give
\[
 \begin{aligned}
 N_\lambda=1-\frac{n-1}{6}L_2\lambda^2+o(\lambda^2),\qquad
 E_\lambda=S_{n,p}^{-p}-\frac{n-1}{6}G_2\lambda^2+o(\lambda^2).
 \end{aligned}
\]
With $t=r^{p/(p-1)}$, the beta recurrence formula gives
\[
 \begin{aligned}
 \frac{S_{n,p}^pG_2}{L_2}
 &=\frac{(n+2)(p-1)}{n-3p+2}\frac{n-p}{n(p-1)}
   =\frac{(n+2)(n-p)}{n(n-3p+2)},\\
 S_{n,p}^pG_2-\frac{n-p}{n}L_2
 &=\frac{3p(n-p)}{n(n-3p+2)}L_2>0.
 \end{aligned}
\]
Hence
\[
 \begin{aligned}
 N_\lambda^{(n-p)/n}-S_{n,p}^pE_\lambda
 &=\frac{n-1}{6}\left(S_{n,p}^pG_2-\frac{n-p}{n}L_2\right)
                  \lambda^2+o(\lambda^2)=c\lambda^2+o(\lambda^2),\quad c>0.
 \end{aligned}
\]
By radial integration,
\[
 \begin{aligned}
 P_\lambda
 &\leq C\lambda^p\left(1+\int_1^{R/\lambda}t^{p-\alpha-1}\,dt\right)=O\bigl(\lambda^{\min\{p,\alpha\}}(1+|\log\lambda|)\bigr)
 =o(\lambda^2).
 \end{aligned}
\]

\emph{The logarithmic case $n=3p-2$.} Here $\alpha=2$ and
$\int_{\Rn}|y|^2V^{\pc}\,dy<\infty$, whereas
\[
 \begin{aligned}
 \int_{|y|<L}|y|^2|\nabla V|^p\,dy
 &=c_0\int_1^L\frac{dt}{t}+O(1)=c_0\log L+O(1),\quad c_0>0.
 \end{aligned}
\]
The cutoff error is $O(\lambda^\alpha)=O(\lambda^2)$. Since $p>2$,
\[
 \begin{aligned}
 C\lambda^4\int_0^{R/\lambda}r^{n+3}|V'(r)|^p\,dr
 &\leq C\lambda^4\left(1+\int_1^{R/\lambda}r\,dr\right)
 =O(\lambda^2),\\
 P_\lambda
 &\leq C\lambda^p\left(1+\int_1^{R/\lambda}r^{p-3}\,dr\right)
 =O(\lambda^2).
 \end{aligned}
\]
Consequently,
\[
 \begin{aligned}
 N_\lambda&=1+O(\lambda^2),\\
 E_\lambda&=S_{n,p}^{-p}
       -\frac{n-1}{6}c_0\lambda^2\log(1/\lambda)+O(\lambda^2),\\
 N_\lambda^{(n-p)/n}-S_{n,p}^pE_\lambda
 &=\frac{n-1}{6}S_{n,p}^pc_0\lambda^2\log(1/\lambda)+O(\lambda^2).
 \end{aligned}
\]
In both cases,
\[
 \frac{N_\lambda^{(n-p)/n}-S_{n,p}^pE_\lambda}
      {2^{p/n}P_\lambda}\longrightarrow+\infty.
\]
By \eqref{eq:counterexample-ratio}, no $B<\infty$ satisfies
\eqref{eq:main-p} when $n\geq3p-2$.
\end{proof}

\bigskip

\noindent \textbf{Acknowledgement.}  R.L. Frank would like to thank T. K\"onig for many extended discussions on the $p=2$ case considered in this paper, which critically shaped his understanding.

\bigskip

\noindent \textbf{AI usage statement.} The authors used ChatGPT 5.6 Sol for exploring proof strategies and developing mathematical arguments, as well as for language editing and the organization of the exposition. The authors take full responsibility for all mathematical statements, arguments, citations, and conclusions in this paper.

\bigskip
\bigskip

\noindent R. L. Frank

\noindent Mathematisches Institut, Ludwig-Maximilians-Universit\"at M\"unchen, Theresienstr. 39, 80333 M\"unchen, Germany, and Department of Mathematics, California Institute of Technology, Pasadena, CA 91125, USA\\
Email: \textsf{frank@math.lmu.de, rlfrank@caltech.edu}

\bigskip

\noindent T. Jin

\noindent Department of Mathematics, The Hong Kong University of Science and Technology\\
Clear Water Bay, Kowloon, Hong Kong\\
Email: \textsf{tianlingjin@ust.hk}

\bigskip

\noindent H.-Y. Wang

\noindent School of Mathematical Sciences, Laboratory of Mathematics and Complex Systems, MOE\\ Beijing Normal University, 
Beijing 100875, China\\
Email: \textsf{wanghuayang@amss.ac.cn}

\bigskip

\noindent J. Xiong

\noindent School of Mathematical Sciences, Laboratory of Mathematics and Complex Systems, MOE\\ Beijing Normal University, 
Beijing 100875, China\\
Email: \textsf{jx@bnu.edu.cn}


\begin{bibdiv}
\begin{biblist}[\small]

\bib{Aubin1976}{article}{
   author={Aubin, Thierry},
   title={Probl\`emes isop\'erim\'etriques et espaces de Sobolev},
   journal={J. Differential Geom.},
   volume={11},
   date={1976},
   number={4},
   pages={573--598},
}

\bib{Aubin1979}{article}{
   author={Aubin, Thierry},
   title={Meilleures constantes dans le th\'eor\`eme d'inclusion
   de Sobolev et un th\'eor\`eme de Fredholm non lin\'eaire pour la
   transformation conforme de la courbure scalaire},
   journal={J. Funct. Anal.},
   volume={32},
   date={1979},
   number={2},
   pages={148--174},
}

\bib{AubinLi1999}{article}{
   author={Aubin, Thierry},
   author={Li, Yan Yan},
   title={On the best Sobolev inequality},
   journal={J. Math. Pures Appl. (9)},
   volume={78},
   date={1999},
   number={4},
   pages={353--387},
}


\bib{BianchiEgnell1991}{article}{
   author={Bianchi, Gabriele},
   author={Egnell, Henrik},
   title={A note on the Sobolev inequality},
   journal={J. Funct. Anal.},
   volume={100},
   date={1991},
   number={1},
   pages={18--24},
}



\bib{ChangYang1991}{article}{
   author={Chang, Sun-Yung A.},
   author={Yang, Paul C.},
   title={A perturbation result in prescribing scalar curvature on $S^n$},
   journal={Duke Math. J.},
   volume={64},
   date={1991},
   number={1},
   pages={27--69},
}

\bib{Dolbeault2021}{article}{
   author={Dolbeault, Jean},
   title={Functional inequalities: nonlinear flows and entropy methods
   as a tool for obtaining sharp and constructive results},
   journal={Milan J. Math.},
   volume={89},
   date={2021},
   number={2},
   pages={355--386},
}



\bib{DolbeaultEstebanLoss2017}{article}{
   author={Dolbeault, Jean},
   author={Esteban, Maria J.},
   author={Loss, Michael},
   title={Interpolation inequalities on the sphere: linear vs. nonlinear flows},
   journal={Ann. Fac. Sci. Toulouse Math. (6)},
   volume={26},
   date={2017},
   number={2},
   pages={351--379},
}


\bib{Druet1998}{article}{
   author={Druet, Olivier},
   title={Optimal Sobolev inequalities of arbitrary order on
   compact Riemannian manifolds},
   journal={J. Funct. Anal.},
   volume={159},
   date={1998},
   number={1},
   pages={217--242},
}

\bib{Druet1999}{article}{
   author={Druet, Olivier},
   title={In\'egalit\'es de Sobolev optimales : r\'ecents
   d\'eveloppements},
   journal={S\'emin. Th\'eor. Spectr. G\'eom.},
   volume={17},
   date={1998--1999},
   pages={111--127},
}


\bib{Druet1999BestConstants}{article}{
   author={Druet, Olivier},
   title={The best constants problem in Sobolev inequalities},
   journal={Math. Ann.},
   volume={314},
   date={1999},
   number={2},
   pages={327--346},
}


\bib{FigalliZhang2022}{article}{
   author={Figalli, Alessio},
   author={Zhang, Yi Ru-Ya},
   title={Sharp gradient stability for the Sobolev inequality},
   journal={Duke Math. J.},
   volume={171},
   date={2022},
   number={12},
   pages={2407--2459},
}

\bib{HangWang2022}{article}{
   author={Hang, Fengbo},
   author={Wang, Xiaodong},
   title={Improved Sobolev inequality under constraints},
   journal={Int. Math. Res. Not. IMRN},
   volume={2022},
   date={2022},
   number={14},
   pages={10822--10857},
}


\bib{HebeyVaugon1996}{article}{
   author={Hebey, Emmanuel},
   author={Vaugon, Michel},
   title={Meilleures constantes dans le th\'eor\`eme d'inclusion de Sobolev},
   journal={Ann. Inst. H. Poincar\'e Anal. Non Lin\'eaire},
   volume={13},
   date={1996},
   number={1},
   pages={57--93},
}


\bib{Li1996}{article}{
   author={Li, Yan Yan},
   title={Prescribing scalar curvature on $\mathbb{S}^n$ and related problems,
   Part II: Existence and compactness},
   journal={Comm. Pure Appl. Math.},
   volume={49},
   date={1996},
   number={6},
   pages={541--597},
}

\bib{Lions1985}{article}{
   author={Lions, Pierre-Louis},
   title={The concentration-compactness principle in the
   calculus of variations. The limit case. I},
   journal={Rev. Mat. Iberoam.},
   volume={1},
   date={1985},
   number={1},
   pages={145--201},
}

\bib{Putterman2020}{article}{
   author={Putterman, Eli},
   title={Cubature formulas and Sobolev inequalities},
   journal={arXiv:2012.08109},
   date={2020},
   url={https://arxiv.org/abs/2012.08109},
}

\bib{SunZou2025}{article}{
   author={Sun, Jingxuan},
   author={Zou, Wenming},
   title={Some improved almost sharp inequalities under constraints},
   journal={Int. Math. Res. Not. IMRN},
   volume={2025},
   date={2025},
   number={16},
   pages={Paper No. rnaf242, 27pp},
}

\bib{Talenti1976}{article}{
   author={Talenti, Giorgio},
   title={Best constant in Sobolev inequality},
   journal={Ann. Mat. Pura Appl. (4)},
   volume={110},
   date={1976},
   pages={353--372},
}

\bib{WYY}{article}{
   author={Wei, Juncheng},
   author={Yan, Zetian},
   author={Ye, Zikai},
   title={On Dolbeault's conjecture},
   journal={arXiv:2610.14359},
   date={2026},
   url={https://arxiv.org/abs/2610.14359},
}

\bib{WangZhang2026}{article}{
   author={Wang, Guofang},
   author={Zhang, Mingwei},
   title={Another sharp criterion for Dirac zero modes},
   journal={arXiv:2609.34958v1},
   date={2026},
   url={https://arxiv.org/abs/2609.34958v1},
}

\end{biblist}
\end{bibdiv}
\end{document}